\documentclass[12pt]{article}
\usepackage[leqno]{amsmath}
\usepackage{amsthm,amscd}
\usepackage{amsfonts} 
\usepackage{amssymb} 

\newtheorem{theorem}{Theorem}[section]
\newtheorem{prop}[theorem]{Proposition}
\newtheorem{corollary}[theorem]{Corollary}

\newtheorem{rem}[theorem]{Remark}
\newtheorem{example}[theorem]{Example}

\def\Proof{\noindent{\sl Proof.}\qquad}

\def\rank{\mathop{\rm rank}}

\newcommand{\mbfz}{{\mathbf z}}

\newcommand{\period}{{\rho}}

\def\kmms{\kern-\mathsurround}

\newcommand{\calM}{{\cal M}}

\newcommand{\calA}{{\cal A}}

\newcommand{\bbR}{{\mathbb R}}

\newcommand{\bbZ}{{\mathbb Z}}

\begin{document}

\title{The characteristic quasi-polynomials of the arrangements of 
root systems
and mid-hyperplane arrangements
}
\author{
Hidehiko Kamiya
\footnote
{
{\it Faculty of Economics, Okayama University}
} 
\\ 
Akimichi Takemura 
\footnote
{
{\it Graduate School of Information Science and Technology,
University of Tokyo} 
} 
\\
Hiroaki Terao 
\footnote
{
{\it Department of Mathematics, Hokkaido University}
}
}
\date{September 2007}
\maketitle

\begin{abstract}
Let $q$ be a positive integer.
In \cite{ktt}, we proved that 
the cardinality of the complement 
of an integral arrangement, after the 
modulo $q$ reduction,
is a quasi-polynomial 
of $q$, which we call the characteristic quasi-polynomial. 
In this paper, we study general properties of the characteristic 
quasi-polynomial as well as discuss two important examples:
the arrangements of reflecting hyperplanes 
arising from irreducible root systems 
and the mid-hyperplane arrangements.
In the root system case,
we present a beautiful formula
for the generating function of the 
characteristic quasi-polynomial which has
been essentially obtained by
Ch.~Athanasiadis \cite{ath96}
and by A.~Blass and B.~Sagan \cite{bls}. 
On the other hand, it is hard to 
find the generating function of the 
characteristic quasi-polynomial in the
mid-hyperplane arrangement case.
We determine them when the dimension is less than six.

\end{abstract}

\section{Introduction}



Let $S$ be an arbitrary $m\times n$ integral matrix 
without
zero columns.
For each positive integer $q\in \bbZ_{>0}$, 
denote $\bbZ_q = \bbZ/q\bbZ$ and 
$\bbZ_{q}^{\times}=\bbZ_q \setminus \{0\}$.
Consider the set  
$$
M_{q} (S) :=
\{
\mbfz = (z_{1}, \ldots, z_{m})\in \bbZ_{q}^{m}: 
\mbfz S \in (\bbZ_{q}^{\times})^{n}
\},$$ 
and its cardinality $|M_{q} (S)|$.
In our recent paper \cite{ktt}, 
we showed that
there exists a monic {\bf quasi-polynomial}
(periodic polynomial) $\chi_{S} (q)$ 
with integral coefficients 
of degree $m$ 
such that
\[
\chi_{S} (q)
=
|M_{q} (S)|, 
\quad q\in \bbZ_{>0}. 
\]
%
Note that the set $M_{q} (S)$ is the complement of an 
arrangement of hyperplanes in the following sense:
Let $S_{1}, S_{2}, \dots, S_{n}$ be the columns of $S$.
Each set
\[
H_{i, q} := \{
\mbfz = (z_{1}, \ldots, z_{m})\in \bbZ_{q}^{m}: 
\mbfz S_{i}  = 0
\},
\quad 1\leq i\leq n, 
\]
can be called a ``hyperplane'' in $\bbZ_{q}^{m}$
by a slight abuse of terminology. Then
$$M_{q} (S) = \bbZ_{q}^{m}  \setminus \bigcup_{i=1}^{n} H_{i, q}.$$  
For a sufficiently large prime number $q$,   
$\chi_{S} (q)$ is known \cite{ath96} to be equal to 
the {\bf characteristic polynomial} \cite[Def. 2.52]{ort} of the
real arrangement consisting of the following hyperplanes
(ignoring possible repetitions):
\[
H_{i, \bbR} := 
\{
\mbfz = (z_{1}, \ldots, z_{m})
\in \bbR^{m}: 
\mbfz S_{i}  
=
0
\},
\quad 1\leq i\leq n. 
\]
It is thus natural to call the quasi-polynomial
 $\chi_{S} (q)$ 
the {\bf characteristic quasi-polynomial} of $S$ 
as in \cite{ktt}.
Let us define 
its 
{\bf
generating function
}
\[
\Phi_{S} (t) :=
\sum_{q = 1}^{\infty} 
\chi_{S} (q) t^{q}.  
\]
We understand that  $M_1(S)=\emptyset$ for $q=1$ and hence 
the summation is in effect for $q\ge 2$.

In this paper, we study the characteristic quasi-polynomial
 $\chi_{S} (q)$ 
or equivalently
its generating function 
$
\Phi_{S} (t).
$
In Section 2, we discuss general properties of the characteristic
quasi-polynomials and their generating functions.  
In the subsequent chapters, we deal with two kinds of specific
arrangements defined over $\bbZ$:
the arrangements of reflecting hyperplanes 
arising from irreducible root systems (Section 3)
and the mid-hyperplane arrangements (Section 4).
Let $R$ be an irreducible root system of rank
 $m$
and $n = |R|/2$.
We assume that an $m\times n$ integral
matrix $S = S(R) = [S_{ij}]$
satisfies
\[
R_{+}  =
\{
\sum_{i=1}^{m}  S_{ij} \alpha_{i} : j = 1, \dots, n
\},
\]
where
$R_{+} $ is a set of positive roots and 
$B(R) = \{\alpha_{1}, \alpha_{2}, \dots, \alpha_{m}\}$ is
the set of {\bf simple roots} associated with $R_{+}$.
In other words, $S$ is a coefficient matrix of
$R_{+} $ with respect to the basis $B(R)$. 
Define 
the characteristic quasi-polynomial
$
\chi_{R} (q) := \chi_{S} (q)
$
and
the generating function
$
\Phi_{R} (t) := \Phi_{S} (t)
$
for each irreducible root system $R$. 
Then $\chi_{R} (q)$ and 
$
\Phi_{R} (t) := \Phi_{S} (t)
$
depend only upon $R$. 
In Section 3, we present a beautiful formula for
the generating function
$\Phi_{R} (t)$
for every irreducible root system $R$. 
This formula has been essentially proved by
Ch.~Athanasiadis in \cite{ath96}
and
A.~Blass and B.~Sagan in \cite{bls}.
(See \cite{hai} also.)
In Theorem 3.1, we will state the formula in our language and
include a proof following
\cite{ath96, bls, hai}
for completeness.  
In Section 4, we will give a formula for 
the generating function
$\Phi_{S} (t)$
when $S$ is equal to the coefficient matrix for
the
mid-hyperplane arrangement
of
dimension less than six.

We are aided by the computer package
PARI/GP \cite{pg}.
\section{Results on the characteristic quasi-polynomial of an integral matrix}

Let $\chi_{S} (t)$
be the characteristic quasi-polynomial of
an $m\times n$ integral matrix $S$
without zero columns.
Fix a nonempty $J \subseteq [n] := \{1,2,\ldots,n\}$ and 
define an $m\times |J|$ matrix $S_{J} $ consisting of 
the columns of $S$ corresponding to the set $J$.  
Let 
$e_{J,1}, \ldots, e_{J,\ell(J)}\in \bbZ_{>0}$
be the elementary divisors of $S_{J}$  
numbered so that 
$e_{J,1} | e_{J,2} | \cdots | e_{J,\ell(J)}$,  
where $\ell(J):=\rank S_J$. 
Write $e(J) := e_{J,\ell(J)}$, and 
define the {\bf lcm period} $\rho_{0} (S)$ of $S$ by
\begin{eqnarray*}
\rho_{0} =
\rho_{0} (S)
&:=& 
{\rm lcm}
\{
e(J) : J \subseteq [n],  
\ J \ne \emptyset 
\} \\ 
&=& 
{\rm lcm}
\{
e(J) : J \subseteq [n],  
\ 1\le |J| \le \min\{m, n\} 
\}. 
\end{eqnarray*}
Then it is known (\cite[Theorem 2.4]{ktt}) that the lcm period
$\rho_{0} $ is a period of 
$\chi_{S} (t)$. 

It is further shown in \cite{ktt} that the constituents of the 
quasi-polynomial $\chi_{S}(t)$ 
are the same 
for all $q$'s with the same value of $\gcd\{ \rho_0, q \}$. 
Let $d$ be a positive integer which divides $\rho_{0}$, and 
define a monic polynomial $P_d(t)=P_{S,d}(t)$ with integral coefficients 
of degree $m$ 
by
\begin{equation}
\label{eq:Pd}
\chi_{S} (q) = P_{d} (q) 
\quad 
\text{for all \ } q \in d + \rho_{0} \bbZ_{\ge 0}.
\end{equation}
Put 
\[
e(J, d) := 
\prod_{j=1}^{\ell(J)}  
\gcd \{ e_{J, j}, d \}.
\]
Then the following formula was essentially proved in 
our previous paper \cite{ktt}.

\begin{theorem}
\label{th:Pd(t)}
For each $d\in \bbZ_{>0}$ with $d|\rho_0$, the polynomial 
$P_d(t)$ is given by  
\[
P_{d} (t) = 
\sum_{J\subseteq [n]} 
(-1)^{|J|}
e(J, d)
t^{m-\ell(J)}, 
\]
where for $J=\emptyset$,  we understand that $\ell(\emptyset)=0$ and that 
$e(\emptyset, d)=1$. 
\end{theorem}

\Proof
Obtained from \cite[(10)]{ktt} and the inclusion-exclusion principle.
\qed

\medskip

\begin{theorem}[\cite{ktt} Theorem 2.5] 
\label{th:P1(t)}
The polynomial $$P_{1} (t)=\sum_{J\subseteq [n]} (-1)^{|J|}t^{m-\ell(J)}$$ 
is equal to the ordinary characteristic polynomial \cite[Def. 2.52]{ort}
of the 
real arrangement consisting of the hyperplanes
(ignoring possible repetitions)
$
H_{1, \bbR},
H_{2, \bbR},
\dots,
H_{n, \bbR}.
$ 
\end{theorem}


\begin{corollary}
\label{co:deg(P-P)<m-s}
Suppose $d, d'\in \bbZ_{>0}$ both divide $\period_0$, 
and assume the following condition holds true for some 
positive integer $s$: 
$
\gcd\{ e(J), d \}
=
\gcd\{ e(J), d' \}
$ 
for all $J\subseteq [n]$ with $|J|\le s$. 
Then
\[
\deg \{ P_{d} (t) - P_{d'} (t)\}  
< m-s. 
\]
In particular, we have
$
\deg \{ P_{d} (t) - P_{1} (t)\}  
< m-s 
$
if 
$
\gcd\{ e(J), d \}
=
1$ 
for all $J\subseteq [n]$ with $|J|\le s$. 
\end{corollary}

\Proof
We apply Theorems \ref{th:Pd(t)}
and
\ref{th:P1(t)}.
It is enough to show 
$
e(J, d) 
= 
e(J, d') 
$ for 
$J \subseteq [n]$ with $\ell(J) \leq s$.
We can choose 
a subset $J'\subseteq J$ such that 
$\ell(J')
=
|J'| = \ell(J) \leq s$.
Then $
\gcd\{ e(J'), d \}
=
\gcd\{ e(J'), d' \}.
$
Since $e(J)|e(J')$ \cite[Lemma 2.3]{ktt}, 
$
\gcd\{ e(J), d \}=
\gcd\{ e(J), d' \}.
$
This shows
$
e(J, d)
=
\prod_{j=1}^{\ell(J)} \gcd\{ e_{J, j},  d \}
=
\prod_{j=1}^{\ell(J)} 
\gcd\{ e_{J, j}, e(J), d \} $\\
$=
\prod_{j=1}^{\ell(J)} 
\gcd\{ e_{J, j}, e(J), d' \} 
=
\prod_{j=1}^{\ell(J)} \gcd\{ e_{J, j}, d' \} 
=
e(J, d').
$
\qed
%
%
\begin{corollary}
\label{co:deg(P+P-P-P)} 
Suppose that
$d\in \bbZ_{>0}$ and $d'\in \bbZ_{>0}$ both divide $\rho_{0}$
and that $\gcd \{d, d'\} = 1$. 
In addition, we assume the following 
condition holds true 
for some positive integer $s$: 
\begin{equation}
\label{eq:=1or=1}
\gcd \{e(J), d\}=1
\text{ \ or \ }
\gcd \{e(J), d'\}=1
\end{equation}
for all $J\subseteq [n]$ with $|J|\le s$. 
Then
\[
\deg 
\{ 
P_{1}(t)+P_{dd'}(t)-P_{d}(t)-P_{d'}(t)
\} 
<m-s. 
\]
\end{corollary}

\Proof
Suppose $J \subseteq [n]$ with $\ell(J) \leq s$.
It is enough to show 
\[
1+e({J, d d'})-e({J, d})-e({J, d'})=0.
\]
We can choose a subset $J' \subseteq J$ 
such that $\ell(J')
=
|J'| = \ell(J)\leq s$.
Then
either
$\gcd \{e(J'), d\}=1
\text{ \ or \ }
\gcd \{e(J'), d'\}=1
$   
by (\ref{eq:=1or=1}).
Since $e(J)|e(J')$, 
$$
\gcd \{e(J), d\}=1
\text{ \ or \ }
\gcd \{e(J), d'\}=1.$$ 
This shows that  
either
$
e(J, d)=1
\text{ \ or \ }
e(J, d') = 1.
$ 
We finally have
\begin{align*} 
0
&=\{ 1 - e(J, d)\}\{1 - e(J, d')\}
=1-e(J, d)-e(J, d')+e(J, d)e(J, d') \\
&=1-e(J, d)-e(J, d')+e(J, d d'). \notag
\end{align*} 
\qed

\begin{corollary}
\label{co:primepowers}
Suppose that
$d\in \bbZ_{>0}$ and $d'\in \bbZ_{>0}$ both divide $\rho_{0}$
and that $\gcd \{d, d'\} = 1$.
If $e(J)$ are prime powers or one 
for all $J$, 
we have $P_{dd'}(t)=P_{d}(t)+P_{d'}(t)-P_{1}(t)$. 
\end{corollary}

\Proof
Easily follows from Corollary \ref{co:deg(P+P-P-P)}. 
\qed

\medskip

For the rest of this section we discuss general properties of
the
generating functions $\Phi_S(t)$ of characteristic quasi-polynomials.
Let $\omega = \exp(2\pi i/\rho_0)$ which is a primitive $\rho_0$'s root
of unity.  By (\ref{eq:Pd}) 
\[
\Phi_S(t)= \sum_{d=1}^{\rho_0} \Phi_{S,d}(t),
\qquad \Phi_{S,d}(t)= \sum_{s=0}^\infty P_d(d+ \rho_0s) t^{d + \rho_0 s}.
\]
Note that $\Phi_S(\omega^k t)=\sum_{d=1}^{\rho_0} \omega^{kd} \Phi_{S,d}(t)$.
Therefore from the orthogonality relations among powers of $\omega$,
i.e.,\ by the
Fourier inversion, 
$\Phi_{S,d}(t)$ for each $d$ can be recovered from $\Phi_S(t)$ by 
\begin{equation}
\label{eq:inversion}
\Phi_{S,d}(t) = \frac{1}{\rho_0}
\sum_{k=1}^{\rho_0} \omega^{-kd} \Phi_S(\omega^k t).
\end{equation}
This relation will be used in Example \ref{ex:BmCmDm} below.

Taking a common denominator, we can express $\Phi_S(t)$ as
a
rational function
\[
\Phi_S(t) = \frac{Q(t)}{(1-t^{\rho_0})^{m+1}},  \qquad \deg Q < (m+1)\rho_0.
\]
In the numerator $Q(t)$ the powers $t^{d+\rho_0s}$, $s=0,1,\dots$, correspond to
$P_d$.  Therefore as in (\ref{eq:inversion}) for each $d$ we can
extract these powers as
\begin{equation}
\label{eq:Qd}
\Phi_{S,d}(t) = \frac{Q_d(t)}{(1-t^{\rho_0})^{m+1}}, \qquad
Q_d(t)=\frac{1}{\rho_{0} }   \sum_{k=1}^{\rho_0} \omega^{-kd} Q(\omega^k t).
\end{equation}
Let 
$$P_{d} (q) = \sum_{k=0}^{m} 
c_{d, k} q^{k}\,\,\,(c_{d, k} \in \bbZ).  $$ 
Then 
\[
Q_{d}(t) 
= (1-t^{\rho_{0} })^{m+1}   
\sum_{s=0}^{\infty} 
P_{d} (d+\rho_{0} s) t^{d+\rho_{0} s}
= (1-t^{\rho_{0} })^{m+1}   
\sum_{k=0}^{m} 
c_{d, k} 
\sum_{s=0}^{\infty} 
(d+\rho_{0} s)^{k}  t^{d+\rho_{0} s}.
\]
Define
polynomials
$q_{d,k}(t)$ by
\begin{equation} 
\label{eq:qdk} 
\sum_{s=0}^\infty (d+\rho_0 s)^k  t^{d+\rho_0 s}
=\frac{q_{d,k}(t)}{(1-t^{\rho_0})^{k+1}}
\qquad
(d=1,\dots, \rho_0).
\end{equation} 
Then we obtain
\begin{equation} 
\label{eq:Qdcdk} 
Q_{d}(t) 
= 
\sum_{k=0}^{m}  
(1-t^{\rho_{0} })^{m-k}   
c_{d, k} q_{d,k}(t).
\end{equation} 
Now we present
the following proposition, to which we give a proof
because we were not able to find an appropriate reference
in literature.

\begin{prop}
\label{prop:fourier}
Define
$q_{d,k}(t)$
by (\ref{eq:qdk}). 
Let $q_{d,k}^{(j)}(1)$ be their $j$-th  derivatives at $t=1$.  Then
\begin{equation}
\label{eq:derivatives}
0\neq 
q_{1,k}^{(j)}(1)=q_{2,k}^{(j)}(1)=\dots=q_{\rho_0,k}^{(j)}(1)
\qquad
(j=0,\dots, k).
\end{equation}
\end{prop}

\medskip
\Proof
For notational simplity  write $q_d(t)=q_{d,k}(t)$ and let
\[
\tilde q_l(t)=\sum_{d=1}^{\rho_0} \omega^{ld}
q_d(t)
\qquad
(l=1,\dots, \rho_0).
\]
Then the inverse Fourier transform is 
\[
q_d(t) = \frac{1}{\rho_0} \sum_{l=1}^{\rho_0} \omega^{-ld}
\tilde q_l(t)
\qquad
(d=1,\dots, \rho_0).
\]
The $j$-th derivative of this at $t=1$  is 
\[
q^{(j)}_d(1)= \frac{1}{\rho_0} \sum_{l=1}^{\rho_0} \omega^{-ld}
\tilde q^{(j)}_l(1).
\]
It follows that $q^{(j)}_d(1)$ does not depend on $d$ if and only if
\begin{equation}
\label{eq:nec-suf}
\tilde q^{(j)}_l(1)=0
\qquad
(l=1,\dots, \rho_0 -1).
\end{equation}
%
By the use of Eulerian numbers $W(k,h)$ (see Chapter III of
\cite{aigner})
we can write
\begin{equation}
\tilde q_l(t)
=\left(\sum_{h=0}^{k-1} W(k,h)
  (\omega^lt)^{k-h}\right)
\big(1 + \omega^l t + \omega^{2l} t^2 + \cdots + \omega^{(\rho_0-1)l}
t^{\rho_0-1}\big)^{k+1}.
\label{eq:tildeQ}
\end{equation}
Note that 
$0 = 1 + \omega^l  + \omega^{2l} + \cdots + \omega^{(\rho_{0} -1)l}$ for
$1 \le l < \rho_0$.  
Therefore differentiating (\ref{eq:tildeQ}) with respect to $t$, 
we have $\tilde q^{(j)}_l(1) = 0$ 
for $1\le l < \rho_0$ and for $0\le j\le k$.
Thus
$q_{1,k}^{(j)}(1)=q_{2,k}^{(j)}(1)=\dots=q_{\rho_0,k}^{(j)}(1)$.
By summing up (\ref{eq:qdk}) we have
\[
\sum_{d=1}^{\rho_{0} } q_{d, k} (t)
=
(1-t^{\rho_{0} } )^{k+1} \sum_{q=1}^{\infty} q^{k} t^{q}.      
\]
Since the Eulerian numbers are  positive integers, it is not hard to see
that
the right hand side is a polynomial of degree 
$\ge k$ with positive integer coefficients.
Thus
$q_{d,k}^{(j)}(1)$
is not zero for $0\le j\le k, \,\, 1\le d \le \rho_{0}$.
\qed

\bigskip
%

\medskip
Proposition \ref{prop:fourier} and (\ref{eq:Qdcdk}) imply that
\[
P_d(t)=P_{d'}(t) \ \Leftrightarrow \ 
c_{d, k}=c_{d', k}\ {\rm for } \ 0\leq k \leq m 
\ \Leftrightarrow \ 
Q_d^{(j)}(1)=Q_{d'}^{(j)}(1)  \ {\rm for } \ 0\le j\le m.
\]
  Furthermore note that
lower order derivatives of $Q_d$ at $t=1$ determine coefficients of
higher degree terms in $P_d(t)$. 
Therefore in terms of the generating function 
the relations in Corollaries \ref{co:deg(P-P)<m-s}
and \ref{co:deg(P+P-P-P)} can be written as follows:
\[
\deg \{ P_{d} (t) - P_{d'} (t)\}  
< m-s 
\   \Leftrightarrow \ 
Q_d^{(j)}(1)=Q_{d'}^{(j)}(1)
\quad (j=0,1,\dots,s),
\]
\begin{align*}
&\deg 
\{ 
P_{1}(t)+P_{dd'}(t)-P_{d}(t)-P_{d'}(t)
\} 
<m-s\\
& \qquad \qquad  \Leftrightarrow \ 
Q_{1}^{(j)}(1)+Q_{dd'}^{(j)}(1)-Q_{d}^{(j)}(1)-Q_{d'}^{(j)}(1)=0
\qquad 
(j=0,1,\dots,s).
\end{align*}

\section{Arrangements of root systems} 
Let $V$ be an $m$-dimensional Euclidean space 
and $E$ be the affine space underlying $V$.
Let $R$ be an irreducible root system in $V$ 
of rank $m$ and $n= |R|/2$.
Suppose that $R_{+}$ is a set of positive roots
and $B = \{\alpha_{1}, \dots , \alpha_{m}  \}$ 
is the set of simple roots associated with $R_{+} $.
Denote the coefficient matrix (with an arbitrary order of columns) of 
the positive roots $R_+$ with respect to $B$  
by 
$S=[S_{ij}]$, which is an $m\times n$ matrix:
\[
R_{+}  =
\{
\sum_{i=1}^{m}  S_{ij} \alpha_{i} : j = 1, \dots, n
\}.
\]
In this section we give an explicit formula for
the generating function $\Phi_{R} (t) := \Phi_{S} (t)$.
The formula was essentially proved by Ch. Athanasiadis 
\cite{ath96} and by A. Blass and B. Sagan \cite{bls}.
(See \cite{hai} also.)
Let $$\tilde{\alpha} = \sum_{i=1}^{m} n_{i} \alpha_{i}$$
be the highest root.  Then it is well-known that
$h := 1 +  \sum_{i=1}^{m} n_{i} $
is equal to the {\bf Coxeter number};
see \cite[Ch. VI, \S 1, 11. Prop. 31]{bou}. 

\begin{theorem}
\label{thm3.1} 
\[
\Phi_{R} (t) =
\frac{(n_{1}\cdots n_{m}) (m!) t^{h}}{(1-t) {\displaystyle \prod_{i=1}^{m}(1-t^{n_{i} } )}}. 
\]
\end{theorem}

Before proving this formula after \cite{ath96, bls, hai},
we introduce basic concepts.  Let $\beta_{1}, \dots, \beta_{m}  $ 
be the basis for $V$ which is dual to the basis $B$:
$(\alpha_{i}, \beta_{j}  ) = \delta_{ij} $.   Define 
a free abelian group
\[
P(R^{\vee} )
=
\bbZ \beta_{1} 
+
\bbZ \beta_{2} 
+
\dots
+
\bbZ \beta_{m} 
\]
of rank $m$. 
Let
\[
H_{i, k} 
=
\{
x\in E
:
(x, \alpha_{i} )
=
k
\}
\]
for 
$1\leq i\leq n$ and $k\in\bbZ$.
Then
$
\calA_{a} 
=
\{
H_{i, k} :
1\leq i\leq n, \, k\in\bbZ
\}
$  
is an arrangement of 
(infinitely many) affine
hyperplanes in $E$.  The reflection
with respect to
$H_{i, k} $ is denoted by
$s_{i, k} $:
\[
s_{i, k} 
(x)
=
x
-
2
\frac{(x, \alpha_{i} )-k}{(\alpha_{i}, \alpha_{i})}  \alpha_{i}  
\qquad (x\in E).
\]
  The 
{\bf
affine Weyl group}
$W_{a}  $ is the group generated by 
$\{s_{i, k} :
1\leq i\leq n, \, k\in\bbZ
\}$.
Each connected component of 
$E \setminus \bigcup \calA_{a} $
is called an {\bf alcove}. The closure of an alcove is
a fundamental domain of the group $W_{a} $ acting on $E$
\cite[Ch. VI, \S 2, 1]{bou}.
Consider a special alcove
\begin{align*}
C 
&=
\{
x\in E
:
0 < (x, \alpha ) < 1
\,\,
(\alpha\in R_{+} )
\}\\
&=
\{
x\in E
:
0 < (x, \alpha_{i} ) < 1
\,\,
(i = 1,\dots, m),
\,\,\,
(x, \tilde{\alpha} ) < 1
\}\\
&=
\{
\sum_{i=1}^{m} c_{i} \beta_{i}  
:
0 < c_{i} 
\,\,
(i = 1,\dots, m),
\,\,\,
\sum_{i=1}^{m} c_{i} n_{i}< 1
\}
\end{align*}
which is an open simplex with vertices
$0$ and  the $\beta_{i}/n_{i}  $ $(1\leq i\leq m)$.
 
\medskip
{\it Proof of Theorem 3.1.}
(\cite{ath96, bls, hai})
Let $A$ be the parallelepiped defined by 
\[
A
=
\{
\sum_{i=1}^{m} c_{i} \beta_{i}  
:
0 < c_{i} < 1
\,\,
(i = 1,\dots, m)
\}.
\]
Then $A$ is, by definition,
a union of alcoves.  The number of 
alcoves in $A$ is equal to
$$
\frac{{\rm vol}(A)}{
{\rm vol}(C)
}
=
(n_{1} \cdots n_{m} ) (m!);$$
see \cite[Ch. VI, \S 2, 4. (5)]{bou}.
Compute
\begin{align*} 
\chi_{R} (t)
&=
|M_{q}(S)|
=
|
\{
{\bf z} = (z_{1}, \dots, z_{m}  ) 
\in
\bbZ_{q}^{m}  
:
{\bf z} S \in (\bbZ^{\times}_{q}  )^{n} 
\}
|\\
&=
|
\{
\gamma 
\in
A\cap
\frac{1}{q}
P(R^{\vee} )
:
\gamma \not\in \bigcup\calA_{a}
\}
|\\
&=
|
(A\cap
\frac{1}{q}
P(R^{\vee} ))
\setminus
\bigcup\calA_{a}
|\\
&=
(n_{1} \cdots n_{m} ) (m!)\, 
|
C\cap
\frac{1}{q}
P(R^{\vee} )
|.
\end{align*} 
On the other hand,
\begin{align*} 
|
C\cap
\frac{1}{q}
P(R^{\vee} )
|
&=
|
\{
\sum_{i=1}^{m} c_{i} \beta_{i}  
:
 c_{i} \in 
\{
\frac{1}{q}, \frac{2}{q}, \dots, \frac{q-1}{q}
\}
\,\,
(i = 1,\dots, m),
\,\,\,
\sum_{i=1}^{m} c_{i} n_{i}< 1
\}
|\\
&=
|
\{
\sum_{i=1}^{m} c_{i} \beta_{i}  
:
 c_{i} \in \bbZ_{>0}
\,\,
(i = 1,\dots, m),
\,\,\,
\sum_{i=1}^{m} c_{i} n_{i}< q
\}
|.
\end{align*} 
This function is known as the Ehrhart quasi-polynomial
of the open simplex bounded by the coordinate hyperplanes and 
the hyperplane $\sum_{i=1}^{m} c_{i} x_{i} = q $;
see \cite[page 235ff]{sta}.
Thus
\begin{align*}
\Phi_{R} (t)
&=
(n_{1} \cdots n_{m})
(m!)
\sum_{q = 1}^{\infty}   \,
|
C\cap
\frac{1}{q}
P(R^{\vee} )|
t^{q}
\\
&=
\frac{(n_{1} \cdots n_{m})
(m!)
t^{n_{1} +\dots+ n_{m} +1}}
{(1-t) {\displaystyle \prod_{i=1}^{m}(1-t^{n_{i} } )}}
=
\frac{(n_{1} \cdots n_{m})
(m!)
t^{h}}
{(1-t) {\displaystyle \prod_{i=1}^{m}(1-t^{n_{i} } )}}. 
\end{align*} 
This completes the proof.
\hspace{108mm}
$\square$ 

\bigskip

\begin{corollary}
The minimum period of the characteristic quasi-polynomial
for an irreducible root system is equal to
${\rm lcm}(n_{1}, \dots , n_{m})$. 
\end{corollary}

\begin{proof}
The assertion holds true because by Proposition 4.4.1 of \cite{sta}
the minimum period of the 
Ehrhart quasi-polynomial
of the open simplex bounded by the coordinate hyperplanes and 
the hyperplane $\sum_{i=1}^{m} n_{i} x_{i} = q $ is equal to
${\rm lcm}(n_{1}, \dots , n_{m})$. 
\end{proof}

The minimum periods for all irreducible root systems are shown in the 
following table:
\bigskip

\begin{tabular}{l|c|c|c} 
root system & $n_{1} , n_{2} , \dots, n_{m}$  
 & $h = 1+n_{1}+n_{2} +\dots+ n_{m}   $ & minimum period\\
\hline
$A_{m} $ & $1,1, \dots,1, 1$ & $m+1$ &$1$\\
$B_{m} $ & $1, 2, 2, \dots, 2$ &$ 2m$& $2$\\
$C_{m} $ & $2, 2, \dots, 2, 1$ & $2m$&$2$\\
$D_{m} $ & $1, 2, 2, \dots, 2, 1, 1$ &$ 2m-2$& $2$\\
$E_{6} $ & $1, 2, 2, 3, 2, 1 $ & $12$ & $6$\\
$E_{7} $ & $2, 2, 3, 4, 3, 2, 1 $ & $18$ & $12$\\
$E_{8} $ & $2, 3, 4, 6, 5, 4, 3, 2 $ & $30$ & $60$\\
$F_{4} $ & $2, 3, 4, 2 $ & $12$ & $12$\\
$G_{2} $ & $2, 3$ & $6$ & $6$\\
\end{tabular} \\

\begin{rem}
With PARI/GP we checked that for 
every irreducible root system the minimum period coincides with the
lcm period $\rho_{0} (S)$.
\end{rem}

\begin{corollary}
\label{th:Coxeter} 
Let $q$ be a positive integer.
For an irreducible root system $R$ with its Coxeter number $h$,
$\chi_{R} (q) > 0 $ if and only if $q \geq h$.  
\end{corollary}

\begin{proof}
The lowest non-zero term of 
$\Phi_{R} (t)$ is equal to
$(n_{1} \cdots n_{m} ) (m!) t^{h}$. 
\end{proof}

\begin{example}
\label{ex:BmCmDm}
$(B_{m}, C_{m}, D_{m})$
By Theorem \ref{thm3.1} and the table above, we have
the generating functions of type $B_{m}, C_{m}, D_{m}$: 
\[
\Phi_{B_{m} } (t)
=
\Phi_{C_{m} } (t)
=
\frac{2^{m-1} (m!) t^{2m} }{(1-t)^{2} (1-t^{2} )^{m-1}},
\,\,\,
\Phi_{D_{m} } (t)
=
\frac{2^{m-3} (m!) t^{2m-2} }{(1-t)^{4} (1-t^{2} )^{m-3}}.  
\]
Thus 
the characteristic quasi-polynomials of $B_{m} $ 
is the same as 
the characteristic quasi-polynomials of $C_{m} $. 
Since the minimum periods of these three root systems
are all equal to two, there exist four polynomials
$P_{1} (q), P_{2} (q), Q_{1}(q), Q_{2} (q) $ satisfying
\[
\Phi_{1} (t)
:=
\Phi_{B_{m} } (t)
=
\Phi_{C_{m} } (t)
=
\sum_{i = 0}^{\infty}  P_{1} (2i+1) t^{2i+1} 
+ 
\sum_{i = 1}^{\infty}  P_{2} (2i) t^{2i}, 
\]
\[
\Phi_{2} (t)
:=
\Phi_{D_{m} } (t)
=
\sum_{i = 0}^{\infty}  Q_{1} (2i+1) t^{2i+1} 
+ 
\sum_{i = 1}^{\infty}  Q_{2} (2i) t^{2i}.
\]
Compute
\begin{align*} 
2 \sum_{i = 1}^{\infty}  P_{2} (2i) t^{2i}
&=
\Phi_{1} (t)
+
\Phi_{1} (-t)
=
\frac{2^{m-1} (m!) t^{2m} }{(1-t)^{2} (1-t^{2} )^{m-1}}
+
\frac{2^{m-1} (m!) t^{2m} }{(1+t)^{2} (1-t^{2} )^{m-1}}\\
&=
\frac{2^{m} (m!) t^{2m}(1+t^{2} ) }{(1-t^{2} )^{m+1}}
\end{align*} 
and
\begin{align*} 
2 \sum_{i = 1}^{\infty}  Q_{1} (2i+1) t^{2i+1}
&=
\Phi_{2} (t)
-
\Phi_{2} (-t)
=
\frac{2^{m-3} (m!) t^{2m-2} }{(1-t)^{4} (1-t^{2} )^{m-3}}
-
\frac{2^{m-3} (m!) t^{2m-2} }{(1+t)^{4} (1-t^{2} )^{m-3}}\\
&=
\frac{2^{m} (m!) t^{2m+1}(1+t^{2}) }{(1-t^{2} )^{m+1}}
=
2 \sum_{i = 1}^{\infty}  P_{2} (2i) t^{2i+1}.
\end{align*} 
This implies
\[
P_{2} (q) = Q_{1} (q-1).
\]
Since $Q_{1} (q)$ is equal to the ordinary
characteristic polynomial of $D_{m}$
by Theorem \ref{th:P1(t)},
we obtain
\[
P_{2} (q) = Q_{1} (q-1)
= (q-2)(q-4)\dots (q-2m+2)(q-m).
\]
Actually we may derive the following characteristic quasi-polynomials
from the generating functions $\Phi_{1} (t)$ and  $\Phi_{2} (t)$ :
\begin{eqnarray*}
\chi_{B_m}(q) 
&=& 
\chi_{C_m}(q) =  
\begin{cases}
(q-1)(q-3)\cdots (q-2m+1) & \text{if $q$ is odd}, \\ 
(q-2)(q-4)\cdots (q-2m+2)(q-m) & \text{if $q$ is even}, 
\end{cases} \\ 
\chi_{D_m}(q) 
&=& 
\begin{cases}
(q-1)(q-3)\cdots (q-2m+3)(q-m+1) & \text{if $q$ is odd}, \\ 
(q-2)(q-4)\cdots (q-2m+4)
\left\{ q^2-2(m-1)q+\frac{m(m-1)}{2} \right\} & \text{if $q$ is even}. 
\end{cases}
\end{eqnarray*}
\end{example}

\noindent
\begin{rem}
We may also prove $\chi_{B_{m}}(2q) = \chi_{D_{m} } (2q-1)$ 
by constructing a one-to-one correspondence
between
$M_{2q} (S(B_{m} ))$ 
and
$M_{2q-1} (S(D_{m} ))$. 
\end{rem}

\begin{example}
Let $R$ be a root system of type $E_{6} $.
We use PLATE V in \cite{bou} 
to get the 
$6\times 36$ matrix $S = S(E_{6})$:
\begin{align*}
S(E_{6} )
&=
\left[
\begin{array}{ccccccccccccccccccccc}
1&0&0&0&0&0&1&0&0&0&0&1&0&0&0&0&1&1&0&0&0\\
0&1&0&0&0&0&0&1&0&0&0&0&1&1&0&0&1&0&1&1&0\\
0&0&1&0&0&0&1&0&1&0&0&1&1&0&1&0&1&1&1&0&1\\
0&0&0&1&0&0&0&1&1&1&0&1&1&1&1&1&1&1&1&1&1\\
0&0&0&0&1&0&0&0&0&1&1&0&0&1&1&1&0&1&1&1&1\\
0&0&0&0&0&1&0&0&0&0&1&0&0&0&0&1&0&0&0&1&1
\end{array}
\right.
\\
&
\left.
\,\,\,
\begin{array}{ccccccccccccccccc}
1&1&0&0&1&1&0&1&1&0&1&1&1&1&1\\
1&0&1&1&1&1&1&1&1&1&1&1&1&1&2\\
1&1&1&1&1&1&1&2&1&1&2&1&2&2&2\\
1&1&1&2&2&1&2&2&2&2&2&2&2&3&3\\
1&1&1&1&1&1&1&1&1&2&1&2&2&2&2\\
0&1&1&0&0&1&1&0&1&1&1&1&1&1&1
\end{array}
\right].
\end{align*} 
Thus $(n_{1}, n_{2}, n_{3},  n_{4}, n_{5}, n_{6})   = 
(1, 2, 2, 3, 2, 1)$ and $h = 12$.
By Theorem \ref{thm3.1}, we have
the generating function of type $E_{6} $: 
\[
\Phi_{E_{6} } (t)
=
\frac{24\cdot(6!) t^{12} }{(1-t)^{3} (1-t^{2} )^{3} (1-t^{3} )}.  
\]
By expanding this formal power series we have
the characteristic quasi-polynomial of $E_{6} $: 
\begin{eqnarray*}
\chi_{E_6}(q)= 
\begin{cases}
q^6 - 36q^5 + 510q^4 - 3600q^3 + 13089q^2 - 22284q +12320 \\ 
\quad = 
(q-1)(q-4)(q-5)(q-7)(q-8)(q-11),  \\ 
\qquad \qquad \gcd\{ 6, q\}=1, \\ 
q^6 - 36q^5 + 510q^4- 3600q^3 + 13224q^2 - 23904q +16640 \\ 
\quad = 
(q-2)(q-4)(q-8)(q-10)(q^2-12q+26), \\ 
\qquad \qquad \gcd\{ 6, q\}=2, \\ 
q^6-36q^5+510q^4-3600q^3+13089q^2-22284q+12960 \\ 
\quad = 
(q-3)(q-9)(q^4-24q^3+195q^2-612q+480), \\ 
\qquad \qquad \gcd\{ 6, q\}=3, \\ 
q^6-36q^5+510q^4-3600q^3+13224q^2-23904q+17280 \\ 
\quad = 
(q-6)^2(q^4-24q^3+186q^2-504q+480), \\ 
\qquad \qquad \gcd\{ 6, q\}=6. 
\end{cases} 
\end{eqnarray*}
We have computed
$
\{e(J): |J|\le 1\}=
\{e(J): |J|\le 2\}=
\{e(J): |J|\le 3\}=\{1\}$,
$
\{e(J): |J|\le 4\}=
\{e(J): |J|\le 5\}=
\{1,2\}$,
$\{e(J): |J|\le 4\}=
\{1,2,3\}$ and the constituents of the quasi-polynomial are
consistent with
Corollaries \ref{co:deg(P-P)<m-s} and 
\ref{co:deg(P+P-P-P)}.
\end{example}

\noindent
{\it Remark.}
R. Suter \cite{sut} gave essentially the same calculation
for every irreducible root system.

\section{Mid-hyperplane arrangement}

The mid-hyperplane arrangement was defined in \cite{kott} to 
find the number of ``ranking patterns'' generated by the unidimensional 
unfolding model in mathematical psychology (\cite{coo}). 
This arrangement is defined as follows. 

Let $m\ge 4$ 
be an integer. 
We define two kinds of hyperplanes as follows: 
\begin{eqnarray*}
H_{ij}&:=& \{ (\alpha_1,\ldots,\alpha_m)\in \bbR^m: \alpha_i=\alpha_j \}, \quad 1\le i<j\le m, \\ 
H_{ijkl}&:=& \{ (\alpha_1,\ldots,\alpha_m)\in \bbR^m: \alpha_i+\alpha_j=\alpha_k+
\alpha_l \}, \quad (i,j,k,l)\in I_4,  
\end{eqnarray*}
where 
$I_4:=\{ (i,j,k,l): 1\le i<j\le m, \ i < k<l\le m, \ 
\text{$j$ is different from $k$ and $l$} \}$. 
Then the mid-hyperplane arrangement $\calM_m$ is defined as  
\[
\calM_m:=\{ H_{ij} \ (1\le i<j\le m), \ H_{ijkl} \ ((i,j,k,l)\in I_4) \}. 
\]

Now, let $T(\calM_m): m\times n$ 
be the coefficient matrix of $\calM_m$, 
where $n=|\calM_m|=\binom{m}{2}+3\binom{m}{4}$.

\subsection{Characteristic quasi-polynomial 
and generating function 
of $\calM_4$} 

When $m=4$, 
we have 
$\{e(J): |J|\le 1\}=\{ 1 \}, \ 
\{e(J): |J|\le 2\}=
\{e(J): |J|\le 3\}= 
\{e(J): |J|\le 4\}=\{ 1,2\}$, 
and thus $\period_0=2$. 
The characteristic quasi-polynomial is 
\allowdisplaybreaks{
\begin{equation*}
\chi_{T(\calM_4)}(q)= 
\begin{cases}
q^4-9q^3+23q^2-15q=q(q-1)(q-3)(q-5) & \text{ if $q$ is odd,} \\ 
q^4-9q^3+26q^2-24q=q(q-2)(q-3)(q-4) & \text{ if $q$ is even.}  
\end{cases} 
\end{equation*}
From this characteristic quasi-polynomial, we obtain 
\[
\Phi_{T(\calM_4)}(t)=
\frac{48 t^6 (t^3 + 5 t^2 + 7 t + 3)}{ (1-t^2)^5}
= \frac{48 t^6 (t+3)}{(1-t)^5 (1+t)^3}. 
\]

\subsection{Characteristic quasi-polynomial 
and generating function 
of $\calM_5$} 

When $m=5$, we have 
$\{e(J): |J|\le 1\}=\{ 1 \}, \ 
\{e(J): |J|\le 2\}=\{ 1, 2\}, \ 
\{e(J): |J|\le 3\}=\{ 1,2,3 \}, \ 
\{e(J): |J|\le 4\}=
\{e(J): |J|\le 5\}=\{ 1,2,3,4,5,6\}$, 
and thus $\period_0=60$. 
With the help of the relations 
\[
P_d(t)+P_{d'}(t)-P_{dd'}(t)-P_1(t)=
\begin{cases}
0 & \text{for} \ (d,d')=(2,5),(3,5),(4,5),(5,6),(5,12), \\ 
-120t & \text{for} \ (d,d')=(2,3),(2,15),(3,4),(3,10),(3,20),(4,15) 
\end{cases}
\]
(see Corollary \ref{co:deg(P+P-P-P)}),  
we obtained the constituents of 
$\chi_{T(\calM_5)}(q)$ as  
\allowdisplaybreaks{
\begin{eqnarray*}
P_{1}(q)&=& 
q^5-25q^4+215q^3-695q^2+504q \\ 
&=&q(q-1)(q-7)(q-8)(q-9), \\ 
P_{2}(q)&=& 
q^5-25q^4+230q^3-920q^2+1104q \\ 
&=&q(q-2)(q^3- 23q^2+ 184q -552), \\ 
P_{3}(q)&=& 
q^5-25q^4+215q^3-735q^2+864q \\ 
&=&q(q-3)(q-9)(q^2- 13 q +32), \\ 
P_{4}(q)&=& 
q^5-25q^4+230q^3-920q^2+1344q \\ 
&=&q (q-4) (q-6) (q-7)(q-8),  \\
P_{5}(q)&=& 
q^5-25q^4+215q^3-695q^2+600q \\ 
&=&q(q-5)(q^3 - 20 q^2+ 115 q-120), \\ 
P_{6}(q)&=& 
q^5-25q^4+230q^3-960q^2+1584q \\ 
&=&q(q-6)(q^3- 19 q^2+ 116 q -264), \\ 
P_{10}(q)&=& 
q^5-25q^4+230q^3-920q^2+1200q \\ 
&=&q(q-10)(q^3- 15 q^2+ 80q -120),  \\ 
P_{12}(q)&=& 
q^5-25q^4+230q^3-960q^2+1824q \\ 
&=&q(q^4- 25 q^3+ 230 q^2  - 960 q +1824),  \\ 
P_{15}(q)&=& 
q^5-25q^4+215q^3-735q^2+960q \\ 
&=&q(q^4- 25 q^3 + 215 q^2 - 735 q +960), \\ 
P_{20}(q)&=& 
q^5-25q^4+230q^3-920q^2+1440q \\ 
&=&q(q^4- 25 q^3 + 230 q^2 - 920 q +1440),  \\ 
P_{30}(q)&=& 
q^5-25q^4+230q^3-960q^2+1680q \\ 
&=&q(q^4- 25 q^3 + 230 q^2 - 960 q +1680),  \\ 
P_{60}(q)&=& 
q^5-25q^4+230q^3-960q^2+1920q \\ 
&=&q(q^4- 25 q^3 + 230 q^2 - 960 q +1920). 
\end{eqnarray*}
}
%
The generationg function is 
\begin{eqnarray*}
\Phi_{T(\calM_5)}(t)
&=&
240
t^{11} (6t^{20} + 40t^{19} + 112t^{18} + 282t^{17} + 511t^{16} \\ 
&& \qquad \qquad 
+917t^{15} + 1301t^{14} + 1818t^{13} + 2163t^{12} + 2493t^{11} \\ 
&& \qquad \qquad \quad 
+ 2479t^{10} +2462t^9 + 2078t^8 + 1734t^7 + 1263t^6 \\ 
&& \qquad \qquad \qquad 
+ 903t^5 + 523t^4 +308t^3 + 137t^2 + 59t + 11) \\ 
&& \quad 
/\{(1-t)^6 (1+t)^4 (1-t+t^2)^2 (1+t+t^2)^3 \\ 
&& \qquad \qquad \qquad \qquad \qquad \qquad 
(1+t^2)^2 (1+t+t^2+t^3+t^4)^2\}. 
\end{eqnarray*}


\noindent
\begin{rem}
For $\calM_6$, we have 
$\{e(J): |J|\le 1\}=\{ 1 \}, \ 
\{e(J): |J|\le 2\}=\{ 1, 2\}, \ 
\{e(J): |J|\le 3\}=\{ 1,2,3 \}, \ 
\{e(J): |J|\le 4\}=\{ 1,2,3,4,5,6\}, \ 
\{e(J): |J|\le 5\}=\{ e(J): |J|\le 6\}=\{ 1,2,3,4,5,6,7,8,9,10,11\}$, 
and thus $\period_0=
27720$. 
It was computationally infeasible 
for us to obtain $\chi_{T(\calM_6)}(q)$, 
so we cannot tell whether this $\period_0$ is the minimum 
period or not. 
\end{rem}

\end{document}